\documentclass{amsart}
\usepackage{amssymb,amscd,amsthm}
\usepackage{latexsym}
\usepackage[dvips]{graphicx}
\date{\today}

\newtheorem{remark}{Remark}
\newtheorem{lemma}{Lemma}
\newtheorem{proposition}{Proposition}
\newtheorem{corollary}{Corollary}
\newtheorem{definition}{Definition}
\begin{document}
\title[Renormalization and $\alpha$-limit set]
{Renormalization and $\alpha$-limit set for expanding Lorenz maps}

\author[Yi-Ming Ding]{Yi-Ming Ding}

\address{Wuhan Institute of Physics and Mathematics, The Chinese Academy of Sciences, Wuhan 430071, P. R. China,}

\email{ ding@wipm.ac.cn}

\thanks{{\it Mathematical classification (2000):}37E05, 37F25, 54H20.}

\maketitle


\medskip

 \medskip

\begin{abstract}
We show that there is a bijection between the renormalizations and
proper completely invariant closed sets of expanding Lorenz map,
which enable us to distinguish periodic and non-periodic
renormalizations. Based on the properties of periodic orbit of
minimal period, the minimal completely invariant closed set is
constructed. Topological characterizations of the renormalizations
and $\alpha$-limit sets are obtained via consecutive
renormalizations. Some properties of periodic renormalizations are
collected in Appendix.

\end{abstract}

\section{Introduction} \ \

Lorenz equations is a system of ordinary differential equations in
$R^3$ which has been enormous influential in Dynamics, providing
inspiration for the definition of a variety of examples including
the geometric models and H\'enon maps \cite{V}. The Lorenz maps we
study are a simplified model for two-dimensional return maps
associated to the flow of the Lorenz equations.

Numerically studies of the Lorenz equations led Lorenz to emphasize
the importance of {\it sensitive dependence of initial
conditions}---an essential factor of unpredictability in many
systems. The simulations for an open neighborhood suggest that
almost all points in phase space approach to a strange
attractor---the Lorenz attractor.  Afraimovic, Bykov and Sil'nikov
\cite{ABS} and Guckenheimer and Willianms \cite{GW} introduced a
geometric model that is an abstraction of the numerically-observed
features possessed by solution to Lorenz equations. Tucker \cite{T1,
T2}  proved the geometric model is valid, so the Lorenz equations
define a geometric Lorenz flow. Luzzatto, Melbourne and Paccaut
showed that such a Lorenz attractor is mixing \cite{LMP}.

A {\em Lorenz map} on $I=[a,\ b]$ is an interval map $f:I \to I$
such that for some
$c\in (a,\ b)$ we have\\
\indent
(i)   $f$ is strictly increasing on $[a,\ c)$ and on $(c,\ b]$;\\
\indent (ii)  $\lim_{x \uparrow c}f(x)=b$, $\lim_{x \downarrow
c}f(x)=a$.

If, in addition, $f$ satisfies the topological expanding condition

(iii) The pre-images set $C=\cup_{n \ge 0}f^{-n}(c)$ of $c$ is dense
in $I$,

then $f$ is said to be an {\em expanding Lorenz map}, \cite{GH,
GS,HS}.

As was mentioned in \cite{LT}, maps with discontinuities are
extremely natural and important, arising for example in billiards or
as return maps for flows with equilibrium points, and very often in
modeling and applications. Lorenz map admits a discontinuity $c$. It
is convenient to leave $f(c)$ undefined, and regard $c$ as two
points, $c+$ and $c-$, $f(c+)=a$ and $f(c-)=b$ from the definition.
Lorenz map plays an important role in the study of the global
dynamics of families of vector fields near homoclinic bifurcations,
see \cite{MPS1,MPS2,Ro,T2,V} and references therein. The expanding
condition follows from \cite{GS,HS,M}, which is weaker than many
other conditions used in \cite{C,GW, Ro} etc.

Renormalization is a central concept in contemporary dynamics. The
idea is to study the small-scale structure of a class of dynamical
systems by means of a renormalization operator $R$ acting on the
systems in this class. This operator is constructed as a rescaled
return map, where the specific definition depends essentially on the
class of systems. The idea of renormalization for Lorenz map was
introduced in studying simplified models of Lorenz attractor,
apparently firstly in Plamer \cite{Pal} and Parry \cite{P4} (cf.
\cite{FL}). The renormalization operator in Lorenz map family, is
the first return map of the original map to a smaller interval
around the discontinuity, rescaled to the original size. Glendinning
and Sparrow \cite{GS} presented a comprehensive study of the
renormalization by investigating the kneading invariants of
expanding Lorenz map.

\begin{definition}\label{def}
A Lorenz map $f:I \to I$ is said to be {\it renormalizable} if there
is a proper subinterval $[u,\ v] \ni c$ and integers $\ell, r>1$
such that the map $g: [u,\ v] \to [u,\ v]$ defined by
\begin{equation} \label{renormalization}
g(x)=\left \{
\begin{array}{ll}
f^{\ell}(x) & x \in [u,\  c), \\
f^{r}(x) & x \in (c,\  v],
\end{array}
\right.
\end{equation}
is itself a Lorenz map on $[u,\ v]$. The interval $[u, \ v]$ is
called the {\it renormalization interval}.

If $f$ is not renormalizable, it is said to be {\it prime}.
\end{definition}

The renormalization map $g$ is the first return map of $f$ on the
renormalization interval $[u,\ v]$(cf. \cite{MM}). Let $f$ be a
renormalizable  Lorenz map. $f$ may have different renormalizations
(cf. \cite{GS, MM}). A renormalization $g=(f^{\ell}, f^{r})$ of $f$
is said to be {\it minimal} if for any other renormalization
$(f^{\ell'}, f^{r'})$ of $f$ we have $ \ell' \ge \ell$ and $r' \ge
r$ (cf.\cite{GS, KM, MM, SS} etc.).

It is not an easy problem to determine wether $f$ is renormalizable
or not. In fact, it is impossible to check if $f$ is prime or not in
finite steps, because $\ell$ and $r$ in (\ref{renormalization}) may
be large.

In this paper we will investigate the renormalization and
$\alpha$-limit set of expanding Lorenz map. The non-expanding case
is more suitable to state in terms of kneading theory, which is
relegated to another paper. The key observation is that one can
renormalize expanding Lorenz map via its {\it proper completely
invariant closed set}, which turns out to be an $\alpha$-limit set
of some periodic points. For given expanding Lorenz map $f$, there
is a one-to-one correspondence between the renormalizations and the
proper completely invariant closed sets of $f$ (Theorem A). Then we
characterize (Theorem B) the renormalizability of $f$ by
constructing the minimal completely invariant closed set $D$, which
is just the $\alpha$-limit set of the periodic orbit with minimal
period. Since the minimal completely invariant closed set
corresponds to the minimal renormalization of $f$, one can define
the renormalization operator $R$ on the space of (expanding) Lorenz
maps: $Rf$ is the minimal renormalization of $f$. Using the
consecutive actions of renormalization operator, we characterize the
$\alpha$-limit set of each point (Theorem C).

The paper is organized as follows. We state our main results in
Section 2, and establish the correspondence between proper
completely invariant closed set and renormalization of expanding
Lorenz map in Section 3. In Section 4, we study the
renormalizability of expanding Lorenz map by constructing the
minimal completely invariant closed set. We characterize the
$\alpha$-limit sets via consecutive renormalizations in Section 5.
At last, we collect some facts about periodic renormalization in
Appendix.

\section{Main results} \ \

For any nonempty open interval $U \subseteq I$, put
\begin{equation}\label{index}
N(U)=\min \left\{n\ge 0: \exists z \in U \ \ {\rm such \ that}\
f^{n}(z)=c \right\}.
\end{equation}
By the definition of $N(U)$, we have $c \in f^{N(U)}(U)$, $N(U) \le
N(V)$ if $ V \subseteq U $, and
\begin{equation}
N(f^i(U))=N(U)-i, \ \ \ \ \ \ \ \ \ i=0,\ 1, \cdots,\ N(U).
\end{equation}
In fact, $N(U)$ is  the maximal integer such that $f^{N(U)}$ is
continuous on $U$. We can regard $N(U)$ as the index of continuity
for the interval $U$. There exists a unique $z\in U$ such that
$f^{N(U)}(z)=c$ because $f^{N(U)-1}$ is continuous and strictly
increasing on $U$. If $f$ is expanding, $N(U)<\infty$ for all open
interval $U$.

$A\subseteq I$, $ A^{\prime }$ represents for the derived set of
$A$, that is, the accumulation point set of $A,A^{\prime \prime
}=(A^{\prime })^{\prime },A^n=(A^{n-1})^{\prime },n=1,2,\cdots.$ $x
\in I$, $orb(x)$ is the orbit with initial value $x$,
$orb(A):=\cup_{x \in A} orb(x)=\cup_{n \ge 0}f^n(A)$.


Recall that a subset $E$ of $I$ is completely invariant under $f$ if
$$f(E)=f^{-1}(E)=E,$$
and it is proper if $E \neq I$.

{\bf Theorem A.}  {\it Let $f$ be an expanding Lorenz map. There is
a one-to-one correspondence between the renormalizations and proper
completely invariant closed sets of $f$. More precisely, suppose $E$
is a proper completely invariant closed set of $f$, put
\begin{equation} \label{periodic pt}
e_{-}=\sup\{x\in E: x<c\},\ \ \ \ \ e_{+}=\inf\{x\in E: x>c\},
\end{equation}
 and
$$\ell=N((e_-,\ c)),\ \ \ \ \ \ \ \ r=N((c,\ e_+)).$$ Then
\begin{equation} \label{periodic pts}
f^{\ell}(e_-)=e_-,\ \ \ \ f^r(e_+)=e_+
\end{equation}
 and the following map
\begin{equation} \label{completely invariant renormalization}
R_Ef(x)=\left \{ \begin{array}{ll}
f^{\ell}(x) & x \in [f^{r}(c+),\   c) \\
f^r(x) & x \in (c, \ f^{\ell}(c-) ]
\end{array}
\right.
\end{equation}
is a renormalization  of $f$.

On the other hand, if  $g$ is a renormalization  of $f$, then there
exists a unique proper completely invariant closed set $B$ such that
$R_Bf=g$. }

\vspace{0.5cm}

A remarkable property of proper completely invariant closed set is
illustrated by (\ref{periodic pts}): the two closest points to $c$,
from the left and right, are periodic. This property is essential
for us to obtain a renormalization. In their  study on the
renormalization theory of expanding Lorenz map via kneading
invariant, Glendinning and Sparrow \cite{GS} proposed a
combinatorial proof for the existence of such two periodic points.

\begin{definition} Suppose $E$ is a proper completely invariant
closed set of  expanding Lorenz map $f$. The renormalization $R_E f$
defined by (\ref{completely invariant renormalization}) in Theorem A
is called the {\it reormalization associated with $E$}. And $E$ is
called the {\it repelling set} associated to the renormalization
$R_E$. The interval $(e_-,\ e_+)$, with endpoints $e_+$ and $e_-$
defined in (\ref{periodic pt}), is called the {\it critical
interval} of $E$ and $R_E$.
\end{definition}

\begin{definition}
A renormalization is said to be periodic if the endpoints of its
{\it critical interval} belong to the same periodic orbit.
\end{definition}

The periodic renormalization is interesting because
$\beta$-transformation
$$T_{\beta, \alpha}(x)=\beta x +\alpha \ \mod
1, \ 1<\beta \le 2, \ 0 \le \alpha<1$$ can only be periodically
renormalized (see \cite{CD2,G},  and Appendix for details). This
kind of renormalization was studied by Alsed$\grave{a}$ and
Falc$\grave{o}$ \cite{AF}, Malkin \cite{M}. It was called phase
locking renormalization by Alsed$\grave{a}$ and Falc$\grave{o}$ in
\cite{AF} because it appears naturally in Lorenz map whose
rotational interval degenerates to a rational point. As we shall see
in Theorem B and Theorem C, the periodic renormalization corresponds
to completely invariant closed set with isolated points, while
non-periodic renormalization corresponds to Cantor set. It is easy
to check if the minimal renormalization is periodic or not (see
Appendix).


By Theorem A, a possible way to characterize the renormalizability
is to look for the {\it minimal completely invariant closed set} $D$
of $f$, in the sense that $D \subseteq E$ for each completely
invariant closed set $E$ of $f$. If we can find a minimal completely
invariant closed set $D$ of $f$, then $f$ is renormalizable if and
only if $D \neq I$. The construction of minimal completely invariant
closed set seems difficult, because we do not even know wether a
Lorenz map always admits such a minimal completely invariant closed
set or not.

The construction of minimal completely invariant closed set is
closely related to the {\it locally eventually onto (l.e.o.)}
property of $f$.

\begin{definition}\label{leo}
Let $f$ be an expanding Lorenz map on $I$. A closed set $B \subseteq
I$ is said to be {\it locally eventually onto (l.e.o.)} under $f$,
if for each open interval $U$ with $U \cap B \neq \emptyset$, there
exists integer $m$ such that $\bigcup_{i=0}^mf^i(U)=I$. And $f$ is
{\it l.e.o.} if $I$ itself is {\it l.e.o.} under $f$.
\end{definition}

We shall construct the minimal completely invariant closed set for
expanding Lorenz map $f$ by choosing some periodic point $p \in I$
and showing that the $\alpha$-limit set of $p$, $\alpha(p)$, is
indeed the minimal completely invariant closed set. The periodic
orbit with minimal period is important in constructing the minimal
completely invariant closed set. It relates naturally to the so
called {\it primary cycle} which was used to characterize the
renormalization of $\beta$-transformation \cite{G}. We begin with
the minimal period $\kappa$ of periodic points of expanding Lorenz
map $f$ (see Lemma \ref{minimal period}). Then we show (see Lemma
\ref{periodic orbit}) that the periodic orbit $O$ with minimal
period is unique, and the $\kappa$-periodic orbit $O$ is {\it
l.e.o.} under $f$. Based on the locally eventually onto property of
$O$, we can prove that the $\alpha$-limit set of each
$\kappa$-periodic point is the minimal completely invariant closed
set.


\vspace{0.5cm}

{\bf Theorem B. } {\it Let $f$ be an expanding Lorenz map with
minimal period $\kappa$, $1<\kappa<\infty$, $O$ be the unique
$\kappa$-periodic orbit, and $D=\overline{\bigcup_{n\ge
0}f^{-n}(O)}$. Then we have the following statements:
\begin{enumerate}
\item $D$ is the minimal completely invariant closed set of $f$.
\item $f$ is renormalizable if and only if $D \neq I$. If $f$ is
renormalizable, then $R_D$, the renormalization associated to $D$,
is the minimal renormalization of $f$.
\item The following trichotomy holds:
\begin{itemize}
\item $ D=I$  if and only if $f$ is prime;

\item $ D=O $  if and only if $R_D$ is periodic;

\item $D$ is a Cantor set if and only if $R_D$ is not
periodic.
\end{itemize}
\end{enumerate}
}

It is easy to see the cases $\kappa=1$ and $\kappa=\infty$ are
prime, Theorem B describes the renormalizability of expanding Lorenz
map completely.

It follows from Theorem B that when $1<\kappa(f)< \infty$, $f$ is
prime if and only if $D=I$. Since $O$ is {\it l.e.o.} under $f$
implies $D=\overline{\bigcup_{n\ge 0}f^{-n}(O)}$ is {\it l.e.o.}
under $f$, $f$ is prime if and only if $f$ is {\it l.e.o.}.

The dynamics of prime expanding Lorenz map $f$ is well understood:
$f$ is prime if and only if it is {\it l.e.o.} (see Corollary
\ref{leo-1}). The {\it l.e.o.} property is an ideal topological
property, which implies transitivity.

Glendinning and Sparrow  \cite{GS} defined the {\it l.e.o.} property
of $f$ in more strict sense. Our definition of {\it l.e.o.} reduces
to their definition when $\kappa \le 2$(see Proposition
\ref{equivalence} in Section 4.3).


According to Theorem B, the minimal renormalizaion of renormalizable
expanding Lorenz map always exists. We can define a renormalization
operator $R$ from the set of renormalizable expanding Lorenz maps to
the set of expanding Lorenz maps (cf. \cite{GH}, \cite{GS}). For
each renormalizable expanding Lorenz map, we define $Rf$ to be the
minimal renormalization map of $f$. For $n>1$, $R^nf=R(R^{n-1}f)$ if
$R^{n-1}f$ is renormalizable. And $f$ is $m$ ($0 \le m \le \infty$)
{\it times renormalizable} if the renormalization process can
proceed $m$ times exactly.  For $0 <i \le m$, $R^if$  is the $i$th
renormalization of $f$. Formally, we denote $R^0f:=f$ as the $0$th
renormalization, whose renormalization interval is denoted by
$[a_0,\ b_0]:=[a,\ b]$.

The consecutive renormalization process can be used to characterize
all the $\alpha$-limit sets and obtain a canonical decomposition of
the nonwandering set of expanding Lorenz map.

Remember that the $\alpha$-limit set $\alpha(x)$ of a point $x \in
I$ under $f$ is defined as
$$
\alpha(x)=\bigcap_{n \ge 0}\overline{\bigcup_{k \ge n}
\{f^{-k}(x)\}}.
$$
$\alpha$-limit set is important in understanding homoclinic behavior
in dynamics. It is often relates to homeomorphsim because the
inverse $\{f^{-1}(x)\}$ is only one point. For endomorphism $f$, the
$\alpha$-limit set is more difficult to understand than
$\omega$-limit set in general, because $\{f^{-k}(x)\}$ is more
complex than $\{f^k(x)\}$. It seems that $\alpha$-limit set is
"difficult" to describe. But $f$ may not have too many different
$\alpha$-limit sets because the $\alpha$-limit set is "large" in
some sense. We have the following unexpected result.

\vspace{0.3cm} {\bf Theorem C.} \label{Theorem C} {\it Let $f$ be an
$m$ ($0\le m \le \infty$) renormalizable expanding Lorenz map,
$[a_i,\ b_i]$ ($0 \le i \le m$) be the  renormalization interval of
the $i$th renormalization $R^if$, and $orb([a_i,\ b_i])=\bigcup_{n
\ge 0}f^n([a_i,\ b_i])$. Then we have:
\begin{enumerate}
\item $f$ admits $m$ proper $\alpha$-limit sets which can be
ordered as
$$\emptyset=E_0 \subset E_1 \subset E_2 \subset \cdot \cdot \cdot \subset E_m \subset I.$$

\item $E_i$ is a Cantor set if the $i$the renormalization is not periodic, and  $E_{i}'=E_{i-1}$ if the $i$the
 renormalization is  periodic.

\item For\ $0<i \le m$, $\alpha(x)=E_i$ if and only if $$x \in orb([a_{i-1},\ b_{i-1}])\backslash orb([a_{i},\
b_{i}]),$$ and $\alpha(x)=I$ if and only if $$x \in
A:=\bigcap_{i=0}^m orb([a_{i},\ b_{i}]).$$
\end{enumerate}
} \vspace{0.3cm}

By Theorem C, we know that expanding Lorenz map admits a cluster of
$\alpha$-limit sets, and we can determine the $\alpha$-limit set of
each point.   Note that $A$ is the attractor of $f$: $A=I$ if
$m=0$,\ $A=orb([a_{m},\ b_{m}])$ if $m < \infty$, and
$A=\bigcap_{i=0}^{\infty}orb([a_i,\ b_i])$ is a Cantor set if
$m=\infty$ (see Theorem D).  So $f$ is prime implies that
$\alpha(x)=I$, $\forall x \in I$. Since $I$ is the largest
$\alpha$-limit set, $f$ admits exactly $m+1$ different
$\alpha$-limit sets.

Remember that the {\it depth} of $A$ is the minimal integer $n$ such
that the $n$th derived set $A^{(n)}$ is empty (cf. \cite{BB}, p.
33). An interesting consequence of Theorem C appears when all the
renormalizations of $f$ are periodic. In this case, Theorem C
implies that, the $i$th derived set of $E_k$ is $E_{k-i}$:
$(E_k)^i=E_{k-i}$ for $0\le i \le k \le m < \infty$. We can
construct closed sets with given depth in a dynamical way.

The proof of Theorem C is based on the 1-1 correspondence between
the $\alpha$-limit sets and completely invariant closed sets: Each
$\alpha$-limit set is a completely invariant closed set (cf. Lemma
\ref{completely invariant}), and each completely invariant closed
set is the $\alpha$-limit set for some periodic point (cf. Lemma
\ref{alpha-limit set}). Using the same ideas, we characterize the
$\alpha$-limit sets of a unimodal map without homterval in
\cite{CD1}.


\vspace{0.3cm}


\vspace{0.5cm}


\section{Completely invariant closed set and renormalization} \ \

Let $f$ be an expanding Lorenz map. A set $E \subseteq I$ is said to
be completely invariant under $f$, if
\begin{equation}\label{completely
invariant0}f(E)=E=f^{-1}(E).\end{equation} A completely invariant
set $E$ is proper if $E \neq I$. Since $f$ is surjective,  $E$ is
completely invariant is equivalent to $E$ is backward and forward
invariant ($ f^{-1}(E)\subseteq E$ and $f(E) \subseteq E$).

Lemma \ref{completely invariant} collects some useful facts of
completely invariant closed set.

\begin{lemma}\label{completely invariant}Let $f$ be an expanding Lorenz
map.
\begin{enumerate}

\item A completely invariant closed set $E$ is proper if and only if
$c \notin  E$;

\item Any proper completely invariant closed set is nowhere dense.

\item The derived set of proper completely invariant closed set is
also completely invariant.

\item $\forall x \in I$, $\alpha(x)$ is a completely invariant closed
set of $f$.

\item  If $p$ is periodic, then $\alpha(p)=\overline{\bigcup_{n\ge
0}f^{-n}(p)}$.

\item If $E$ is a completely invariant closed set of $f$, then for $A \subset I$, we have
\begin{equation}
f^{-1}(A \cap E)=f^{-1}(A) \cap E, \ \ \ \ \ \ \ \ \ f(A \cap
E)=f(A) \cap E.
\end{equation}
\end{enumerate}

\end{lemma}

\begin{proof}
1, It is necessary to prove $c \in E$ implies that $E=I$. By the
invariance of $E$ under $f^{-1}$, $c \in E$ implies that $f^{-n}(c)
\in E$. So $\bigcup_{n \ge 0}f^{-n}(c) \subset E$, which implies
that $E \supseteq \overline{\bigcup_{n \ge 0}f^{-n}(c)}=I$.

2, If $E$ contains some interval $U$, then $c \in f^{N(U)}(U)
\subseteq E$ because $E$ is invariant under $f$, we obtain a
contradiction. So $E$ contains no interval.

3, Suppose $E$ is a proper completely invariant closed set. It
follows that $c \notin E$. So both $f$ and $f^{-1}$ are continuous
at each point of $x \in E$, which implies that $E'$ is backward
invariant and forward invariant.

4, $x \in I$,  $\alpha(x)=\bigcap_{n \ge 0}\overline{\bigcup_{k \ge
n} \{f^{-k}(x)\}}$. For each $n \in N$, $\bigcup_{k \ge n}
\{f^{-k}(x)\}$ is invariant under $f^{-1}$, it follows
$$f^{-1}(\alpha(x)) \subseteq \alpha(x).$$

Remember $y\in \alpha(x)$ is equivalent to the fact that there
exists a sequence $\{x_k\} \subset I$ and an increasing sequence
$\{n_k\} \subset N$ such that $f^{n_k}(x_k)=x$ and $x_k \to y$ as $k
\to \infty$.  Assume $y \in \alpha(x)$,  we have $f(y) \in
\alpha(x)$ if $y$ is not the discontinuity $c$. If $y=c$ we consider
$c$ as two points $c+$ and $c-$. It is easy to see $c+\in \alpha(x)$
implies $f(c+)=0\in \alpha(x)$, and $c-\in \alpha(x)$ implies
$f(c-)=1\in \alpha(x)$. So we conclude
$$f(\alpha(x))\subseteq \alpha(x).$$

5, If $p$ is periodic with period $m$, then $p\in f^{-km}(p)$ for
all $k \in N$, which implies that $p \in \alpha(p)$. Since
$\alpha(p)$ is completely invariant, we know that $f^{-n}(p) \subset
\alpha(p)$. We have $\alpha(p)\supseteq \overline{\bigcup_{n\ge
0}f^{-n}(p)}$. The converse inclusion $\alpha(p) \subseteq
\overline{\bigcup_{n\ge 0}f^{-n}(p)}$ is trivial.

6, Since $E$ is completely invariant, it follows  $f^{-1}(A \cap
E)=f^{-1}(A) \cap f^{-1}(E)=f^{-1}(A) \cap E$. The first equality
holds.

The inclusion  $f(A\cap E) \subseteq f(A) \cap E$ is trivial. To
prove the converse inclusion, suppose $x\in f(A) \cap E$. $x \in
f(A)$ implies that $f(y)=x$ for some $y \in A$. $x \in E$ implies
$\{f^{-1}(x)\} \subseteq E$. As a result, one gets $y \in E$. Hence,
$y \in A\cap E$, which implies $f(A) \cap E \subseteq f(A\cap E)$.
The second equality follows.
\end{proof}

For expanding Lorenz map $f$, Lemma \ref{completely invariant}
indicates that each completely invariant closed set containing $c$
is trivial. It is possible that all the completely invariant closed
set of $f$ is trivial. If this is the case, $f$ is prime because
$\alpha(x)=I$ for all $x \in I$.

\begin{lemma}\label{complimentary}Let $f$ be an expanding Lorenz
map, $E$ be a proper completely invariant closed set of $f$,
$J_E=(e_-,\ e_+)$ be the critical interval of $E$. $N((e_-,\
c))=\ell$, $N((c,\ e_+))=r$, $[u,\ v]=[f^r(c+),\ f^{\ell}(c-)]$.
Then

\begin{equation} \label{comp}
I\backslash E= \bigcup_{n \ge 0}f^{-n}(J_E)= \bigcup_{n \ge
0}f^{-n}((u,\ v)).
\end{equation}
\end{lemma}

\begin{proof} Since $E$ is completely invariant, we have $E \cap \bigcup_{n \ge
0}f^{-n}(J_E)=\emptyset$, which indicates $\bigcup_{n \ge
0}f^{-n}(J_E) \subseteq I\backslash E $.

$x \in I\backslash E$, there exists an open interval $U$ such that
$x \in U \subset I\backslash E$ because $I\backslash E$ is open.
Furthermore, we can assume that $U$ is the maximal open interval
containing $x$ which belongs to $I\backslash E$. Since $f$ is
expanding, $N(U)<\infty$, and $c \in f^{N(U)}(U)$. It follows that
$f^{N(U)}(U) \subset J_E$ because $f^{N(U)}(U)\cap E =\emptyset$.
The maximality of $U$ indicates that $f^{N(U)}(U)=J_E$. As a result,
$f^{N(U)}(x)\in J_E$, i.e., $x \in f^{-N(U)}(J_E)$. Hence,
$I\backslash E \subseteq \bigcup_{n \ge 0}f^{-n}(J_E)$. We have
proved $I\backslash E= \bigcup_{n \ge 0}f^{-n}(J_E)$.

Since $E$ is a completely invariant closed set, we have $(u,\ v)
\subseteq J_E$. It follows that $\bigcup_{n \ge 0}f^{-n}(J_E)
\supseteq  \bigcup_{n \ge 0}f^{-n}((u,\ v))$.

$\forall x \in (e_-,\ c)$, put $\ell_x=N((e_-,\ x))$. We get
$f^{\ell_x}(x) \in (c,\ v)$ by the complete invariance of $E$. So we
conclude $(e_-,\ c) \subset \bigcup_{n \ge 0}f^{-n}((u,\ v))$. By
the same argument, we can obtain $(c,\ e_+) \subset \bigcup_{n \ge
0}f^{-n}((u,\ v))$. So $J_E=(e_-,\ e_+) \subset \bigcup_{n \ge
0}f^{-n}((u,\ v))$. Since $\bigcup_{n \ge 0}f^{-n}((u,\ v))$ is
backward invariant, we have $$\bigcup_{n \ge 0}f^{-n}(J_E) \subseteq
\bigcup_{n \ge 0}f^{-n}((u,\ v)).$$ This completes the proof of
(\ref{comp}). $\hfill \Box$


\subsection{Proof of Theorem A}

Suppose $E$ is a proper completely invariant closed set of $f$.
$e_+$, $e_-$, $\ell$ and $r$ are defined as in the statement of
Theorem A. By Lemma \ref{completely invariant}, $e_-,\ e_+,\ \ell$
and $r$ are well defined and $e_-<c<e_+$.

At first, we prove $f^{\ell}(e_-)=e_-$.

By the definition of $\ell$, $f^{\ell}$ is continuous and monotone
on $(e_-,\ c)$. Put $z$ be the unique point in $(e_-,\ c)$ such that
$f^{\ell}(z)=c$. Since $E$ is completely invariant, we conclude that
$f^{\ell}(e_-)=e_-$. In fact, if $f^{\ell}(e_-)>e_-$, then
$e_-<f^{\ell}(e_-)<f^{\ell}(z)=c$, which contradicts to the
maximality of $e_-$ because $f^{\ell}(e_-) \in E \bigcap (e_-,\ c)$.
If $f^{\ell}(e_-)<e_-$, there must be some point  $y \in (e_-,\ c)$
such that $f^{\ell}(y)=e_-$, which contradicts also to the
maximality of $e_-$ and the complete invariance of $E$ under $f$.

Similarly, we can prove $f^r(e_+)=e_+$.

Since $E$ is completely invariant, we  conclude
$$f^{\ell}((e_-,\ c))=(e_-,\ f^{\ell}(c-)) \subseteq (e_-,\ e_+).$$
If, on the contrary, $f^{\ell}(c-)>e_+$, then there exists $z\in
(e_-,\ c)$ such that $f^{\ell}(z) =e_+$, which implies $z\in E$
because $E$ is completely invariant. This contradicts to the
minimality of $e_+$.

Similarly, $$f^{r}((c,\ e_+)) \subseteq (e_-,\ e_+).$$

It follows that the map $R_Ef$ defined in Theorem A is a
renormalization of $f$.

Now we prove the second statement. Suppose $g=(f^{m}, \ f^k)$ is a
renormalization map of $f$ with renormalization interval $[u,\
v]:=[f^k(c+),\ f^{m}(c-)]$. Put
$$F_g=\{x \in I, orb(x)\cap (u,\ v)\neq \emptyset \},$$
$$J_g=\{x \in I, orb(x)\cap (u,\ v)=\emptyset \}.$$
Since $F_g=\bigcup_{n \ge 0}f^{-n}((u,\ v))$, $F_g$ is a completely
invariant open set. And $J_g=I \backslash F_g$ is a completely
invariant closed set of $f$. $R_{J_g}=g$ follows from Lemma
\ref{complimentary}.

The proof of Theorem A is complete.
\end{proof}

\section{Minimal completely invariant closed set} \ \

Applying Theorem A, the renormalizability problem of expanding
Lorenz map reduces to check wether it admits a proper completely
invariant closed set. In this section, we shall construct the
minimal completely invariant closed set of $f$. We begin with the
minimal period of the periodic orbits of $f$, and show that the
periodic orbit $O$ with minimal period of $f$ is unique. Then we
conclude that periodic orbit $O$ has the {\it locally eventually
onto ( l.e.o.)} property, which enables us to show that the
$\alpha$-limit set $D:=\alpha(O)=\overline{\bigcup_{n \ge
0}f^{-n}(O)}$ is the minimal completely invariant closed set of $f$.
By Theorem A, $f$ is renormalizable if and only if $D \neq I$. Based
on the structure of $D$, we can prove Theorem B. Using Theorem B, we
can obtain two Propositions about the {\it l.e.o.} property and {\it
trivial renormalization} of $f$.

\subsection{Periodic orbit with minimal period} \ \ \ \

In this subsection, we will show that the periodic orbit with
minimal period is very special because it relatives to the minimal
completely invariant closed set.

The period of periodic points of Lorenz map was well studied by
Alsed$\grave{a}$ et al in \cite{ALMT}. It is shown that a Lorenz map
is asymptotically periodic if and only if the derived set $C'(f)$ of
$C(f)=\bigcup_{n \ge 0} f^{-n}(c)$ is countable \cite{DF}. The
following Lemma \ref{minimal period} determines the minimal period
$\kappa$ of periodic points of expanding $f$ via the preimages of
$c$.

\begin{lemma} \label{minimal period} Suppose $f$ is an expanding Lorenz map on $[a,\ b]$ without fixed point.
The minimal period of $f$ is equal to $\kappa=m+2$, where
\begin{equation} \label{mm}
 m=\min\{i \ge 0: f^{-i}(c) \in [f(a), f(b)]\}.
\end{equation}
\end{lemma}

\begin{proof}
We prove the result by two steps:  we first prove that $f$ has
$(m+2)$-periodic point, then we show that $f$ has no periodic point
with period less than $m+2$.

Notice that $f^{-m}(c) \in [f(a), f(b)]$ and $x$ admits two
preimages if and only if $x \in [f(a),\ f(b)]$. Let $c_{m+1}$ and $
c_{m+1}'$ with $c_{m+1} < \ c_{m+1}'$ be the two preimages of
$f^{-m}(c)$. The set $f^{-i}(c)$ for $i=0, 1, \cdots, m$ is a
singleton. Denote $c_i:=f^{-i}(c)$, $i=0, 1, \cdots, m$. Let $Q_1\in
(a,\ c)$ and $Q_2 \in (c,\ b)$ be the points such that $f(Q_1)=f(b)$
and $f(Q_2)=f(a)$. See Figure 1 for an intuitive picture of $m=2$.

\begin{figure}
\centering
\includegraphics[width=0.8\textwidth]{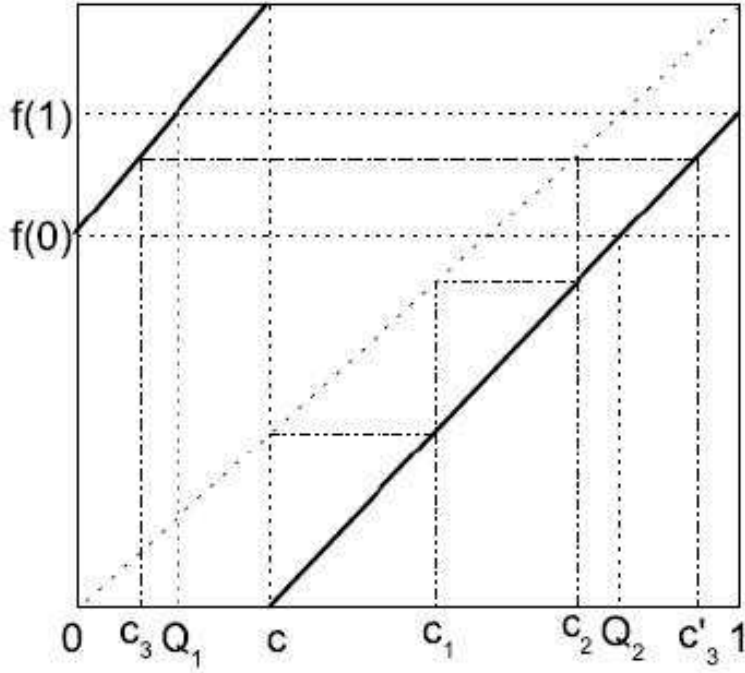}
\caption{A Lorenz map with m=2}
\end{figure}

Since $m$ is the smallest integer such that $f^{-m}(c) \in [f(a),
f(b)]$, we have
\begin{equation} \label{order}
 c_{m+1}\le Q_1 < c_i <  Q_2 \le c_{m+1}' \qquad (0\le i \le m).
\end{equation}

Let $c_{i_0}$ be the minimal point in $\{c, c_1, \cdots, c_m\}$. For
interval $[c_{m+1}, c_{i_0}]$, by (\ref{order}), we obtain that
\begin{eqnarray*}
& & [c_{m+1},\ \ c_{i_0}]\stackrel{f^{i_0}}{\longrightarrow}
[c_{m+1-{i_0}},\ \ c] \stackrel{f}{\longrightarrow} [c_{m-{i_0}},\ \
b] \supseteq [c_{m-{i_0}},\ \
c_{m+1}'] \\
& & [c_{m-{i_0}},\ \
c_{m+1}']\stackrel{f^{m-{i_0}}}{\longrightarrow}[c,\ \
c_{{i_0}+1}]\stackrel{f}{\longrightarrow}[a,\ \ c_{i_0}] \supseteq
[c_{m+1},\ \ c_{i_0}],
\end{eqnarray*}
which implies that
$$
[c_{m+1},\ c_{i_0}] \subseteq f^{m+2}([c_{m+1},\ c_{i_0}]).
$$
So, $f$ has an $m+2$-periodic point in $[c_{m+1}, c_{i_0}]$.

Fix $1<j<m+2$. We shall prove that $f$ admits no $j-$periodic point.
Put $c_{\ell}=\min \{c,\ c_1, \ldots,\ c_m\}$, $c_{r}=\max \{c,\
c_1, \ldots,\ c_m\}$.

{\bf Claim: }  $f$ can not have $j$-periodic points in $(a,\
c_{\ell})$ and $(c_{r},\ b)$.

By the selection of $m$, we get $N((a,\ c_{\ell}))>m$. So $f^j$ is
continuous and monotone on $(a,\ c_{\ell})$. It is easy to see
$f^j(a)>a$. If $f^j$ admits a fixed point $x_*$ in $(a,\ c_{\ell})$,
then
$$a<f^j(a)<f^{2j}(a)< \cdots< f^{nj}(a)<x_*,  \ \ \ \ n>0.$$
So $\{f^{nj}(a)\}_n$ approaches to a fixed point of $f^j$ as $n \to
\infty$, which is impossible because expanding Lorenz map does not
admits attractive periodic orbit.

Similarly, if $f^j$ admits a fixed point in $(c_{r},\ b)$, then
$\{f^{nj}(b)\}_n$ will converge to a fixed point of $f^j$, which
contradicts to $f$ is expanding.

Now, for any open interval $J$ with both endpoints in $\{c,\ c_1,
\ldots,\ c_m\}$ and $J\cap \{c,\ c_1, \ldots,\ c_m\}=\emptyset$,
 we know that $N(J)>m$, and at least one of the following cases hold:
\begin{itemize}
\item $f^j(J) \cap J=\emptyset$;
\item $f^i(J) \subseteq ((a,\ c_{\ell}) \cup (c_{r},\ b))$ for some
$1<i \le j$.
\end{itemize}
It follows that $f$ admits no $j$-periodic point in $J$.
\end{proof}

\begin{remark}
For $m$ defined in (\ref{mm}), it is interesting to note when
$f^{-m}(c)$ is happen to be one of the endpoints of $[f(a), \
f(b)]$. If $f^{-m}(c)=f(a)$, then $c+$, as well as $a$, is a
periodic point with period $m+2$. If $f^{-m}(c)=f(b)$, then
$f^{m+2}(c-)=c-$.



\end{remark}

Let $P_L$ be the largest $\kappa-$periodic point in $[a,\ c)$ and
$P_R$ be the smallest $\kappa-$periodic point in $(c,\ b]$.

\begin{lemma}\label{equalities} \ \ Put $L_1=(P_L,\ c)$, $R_1=(c, P_R)$. We have
\begin{equation}N(L_1)=N(R_1)=\kappa. \end{equation}
\end{lemma}

\begin{proof}
We only prove $N(L_1)=\kappa$ by showing both $N(L_1)<\kappa$ and
$N(L_1)> \kappa$ are impossible.

Suppose that $N(L_1)<\kappa$. We have $f^{N(L_1)}(P_L)\not=P_L$. By
the definition of $N(L_1)$, there exists $z \in L_1$ such that
$f^{N(L_1)}(z)=c$. Since $P_L$ is the largest $\kappa$-periodic
point of $f$ in $[0,\ c)$ and $f^{N(L_1)}(P_L)$ is a
$\kappa$-periodic point,
 we must have
$$
f^{N(L_1)}(P_L)<P_L.
$$
$N(L_1)<\kappa$ implies that $f^{N(L_1)}$ is increasing on
$[P_L,c)$. For the interval $(P_L, \ z)$, it follows that
$$
f^{N(L_1)}([P_L,z))=[f^{N(L_1)}(P_L),c) \supseteq [P_L,z).
$$
So there exists $P_* \in (P_L,z)$ such that $f^{N(L_1)}(P_*)=P_*$ by
the continuity of $f^{N(L_1)}$ on $(P_L,z)$. Hence $P_*$ is a
periodic point of $f$ with period $N(L_1)<\kappa$, which contradicts
to the minimality  of $\kappa$.

Assume $N(L_1)>\kappa$. It follows from (\ref{index}) that
$f^{\kappa}$ is continuous and increasing on $L_1=[P_L,c)$. We have
to exclude two cases: $f^{\kappa}(c-)>c$ and $f^{\kappa}(c-)<c$,
which imply that $N(L_1)>\kappa$ is also impossible.

If $f^{\kappa}(c-)>c$, there exists $z \in (P_L, c)= L_1$ such
$f^{\kappa}(z)=c$, which contradicts to the minimality  of $N(L_1)$.

If $f^{\kappa}(c-)<c$, by the monotone property of $f^{\kappa}$ on
$[P_L, c)$, we obtain a decreasing sequence $\{f^{n \kappa}(c-)\}$
with lower bound $P_L$. Hence,
$$
f^{-n}(c) \cap [P_L, c)= \emptyset,
$$
which contradicts to the fact that $f$ is expanding.
\end{proof}

\begin{lemma}\label{periodic orbit} \ \ Suppose that $f$ is an expanding
 Lorenz map, and
$1<\kappa <\infty$ is the smallest period of the periodic points of
$f$. Then

{\rm i)} \ \  $f$ admits a unique $\kappa-$periodic orbit;

{\rm ii)} We have $$\bigcup_{i=0}^{\kappa-1}f^i([P_L,P_R])=I;$$

{\rm iii)} For any open interval $U$ containing a $\kappa$ periodic
point, there exists positive integer $n$ such that
$$\bigcup_{i =0}^n f^i(U)=I.$$
\end{lemma}

\begin{proof}
i)\ Suppose that $f$ has two distinct $\kappa$-periodic orbits
$orb((P_L)$ and $orb(Q_L)$, where $P_L$ and $Q_L$ are the maximal
points in $L$ of these two periodic orbits respectively. Without
loss of generality, we can suppose $P_L$ is the largest
$\kappa$-periodic point in $L$. Put $L_1=(P_L,\ c)$ and $L_2=(Q_L,\
c)$.

By Lemma \ref{equalities} and $L_1 \subset L_2$ we know that $N(L_2)
\le N(L_1)=\kappa$. If $N(L_2) =\kappa_1<\kappa$,  there exists a
point $z \in (Q_L, P_L)$ such that $f^{\kappa_1}(z)=c$. Since
$f^{\kappa_1}(Q_L)<c$, it follows that $f^{\kappa_1}(Q_L)<Q_L$
according to the choice of $Q_L$. So we have
$$f^{\kappa_1}((Q_L, z))=(f^{\kappa_1}(Q_L), c) \supset (Q_L, \ z),$$ which
implies that $f$ admits an ${\kappa_1}-$periodic point in $(Q_L,
P_L)$. We obtain a contradiction because $\kappa$ is the minimal
period of periodic points. So we conclude that $N(L_2)=\kappa$ and
$f^{\kappa}$ is continuous on $L_2$.

Consider the action of $f^{\kappa}$ on the interval $[Q_L,\ P_L]$,
we have
$$f^{n \kappa}([Q_L,\ P_L])=[Q_L,\ P_L],$$ which contradicts to the
fact $f$ is expanding.

So $f$ admits a unique $\kappa$-periodic orbit.

ii) By the proof of i), we get $f^{\kappa}([P_L,\ c)) \supset [P_L,\
c)$ and $f^{\kappa}((c,\ P_R]) \supset (c,\ P_R]$.  Observe that
$$f([P_L,\ c))=[f(P_L),\ b), \ \ \ \ \ \ \ \
f^2([P_L,\ c))=[f^2(P_L),\ f(b))$$
$$f((c,\ P_R])=(a,\ f(P_R)], \ \ \ \ \ \ \ \
f^2((c,\ P_R])=(f(a),\ f^2(P_R)].$$ Since $f$ is expanding implies
 $f(a) \le f(b)$, we conclude
 $$
 f^2([P_L,P_R]) \supseteq
[f^2(P_L),f^2(P_R)].
$$ It follows
\begin{equation} \label{period}
f^i([P_L,\ P_R])\supseteq [f^i(P_L),\ f^i(P_R)] \ \ \ \ \ for \ \ \
\ i=2,\cdots,\kappa.
\end{equation}

Hence,
$$\bigcup_{i=0}^{\kappa-1}f^i([P_L,\ P_R])=I.$$

iii) Put $U=(x,\ y)$.Without loss of generality, we only consider
the case $P_L \in U$ because some iterates of $U$ contains $P_L$.
Let $N((P_L,\ y))=i$, $N((x,\ P_L))=j$. Since
$$f^{i}([P_L,\ y)) \supseteq [P_L,\ c),\ \ \ \ \ \ \ f^{j}((x,\ P_L]) \supseteq
(c,\ P_R].$$ The conclusion follows from ii).
\end{proof}

The third statement in the Lemma \ref{periodic orbit} implies that
the $\kappa$-periodic orbit $O$ is {\it l.e.o.} under $f$. This
result holds trivially for $D=\overline{\bigcup_{n\ge 0}f^{-n}(O)}$.
This is why the $\alpha$-limit set $D$ of the $\kappa$-periodic
orbit is so important in describing the renormalization of expanding
Lorenz map.

\vspace{0.3cm}




\vspace{0.2cm}

\subsection{Proof of Theorem B} \ \
According to Lemma \ref{periodic orbit} $f$ admits a unique periodic
orbit $O$ with period $\kappa$.   We denote
$D:=\alpha(O)=\overline{\cup_{n \ge 0}f^{-n}(O)}$ as the
$\alpha$-limit set of the $\kappa$-periodic orbit of $f$.

(1) By Lemma  \ref{completely invariant} we know that $D$ is a
completely invariant closed set. We shall prove that $D$ is minimal.

Suppose $E$ is a completely invariant closed set. We have two cases:

{\bf Case 1:} $E\cap (P_L,\ P_R) \neq \emptyset$.

In this case, we can suppose that $(P_L,\ c)\cap E \neq \emptyset$
without loss of generality. Assume that $y \in (P_L,\ c)\cap E$. By
Lemma \ref{equalities} we know that $f^{\kappa}$ is continuous on
$(P_L,\ c)$ and $f^{\kappa}((P_L,\ c))\supset (P_L,\ c)$.  So there
exists $y_1 \in (P_L,\ c)\cap E$ and $y_1<y$ such that
$f^{\kappa}(y_1)=y$. Similarly, we can obtain a decreasing sequence
$\{y_n\}\subset (P_L,\ c)\cap E$ such that
$f^{\kappa}(y_{n+1})=y_n$, $n=1,\ 2, \cdots$ and $\lim_{n \to
\infty}y_n=P_L$. So $P_L \in E$ because $E$ is closed. Hence,
$\cup_{n \ge 0}f^{-n}(O) \subset E$ because $E$ is backward
invariant. $E$ is closed implies that $D \subseteq E$.

{\bf Case 2:} $E\cap (P_L,\ P_R) = \emptyset$.

By Lemma \ref{periodic orbit} \ ii) we know that
$\bigcup_{i=0}^{\kappa-1}f^i([P_L,\ P_R])=I$. So $[P_L,\ P_R] \cap E
\neq \emptyset$. The assumption $E\cap (P_L,\ P_R) = \emptyset$
indicates $[P_L,\ P_R]\cap E=\{P_L,\ P_R\}$, which implies that
$D=E=O$.

The proof of the minimality of $D$ is complete.

(2) By Theorem A, $f$ is renormalizable if and only if $f$ admits a
proper completely invariant closed set. Since $D$ is the minimal
completely invariant closed set, we know that $f$ is renormalizable
is equivalent to $D\neq I$.

If $D\neq I$, according to Theorem A, we know that $R_D$ is a
renormalization of $f$, where
\begin{eqnarray*} \label{renormal}R_D f(x)=\left \{ \begin{array}{ll}
f^{\ell}(x) & x \in [f^{r}(c+),  c) \\
f^{r}(x) & x \in (c,  f^{\ell}(c-) ],
\end{array}
\right.
\end{eqnarray*}
and
\begin{eqnarray*}
&\ell=N([d_-, c))& \ \ \ \ \ \ \  d_-= \sup \{x<c: x \in
D\}, \\
&r=N((c, d_+))&\ \ \ \ \ \ \ d_+= \inf \{x>c: x \in D\}.
\end{eqnarray*}

Assume $g=(f^{\ell'},\ f^{r'})$ is a renormalization of $f$ with
renormalization interval $[u,\ v]$. By Theorem A there exists a
completely invariant closed set
$$
E=\{x \in I: orb(x) \cap (u,\ v)=\emptyset \}
$$
such that $g=R_E$, and
\begin{eqnarray*} &\ell'=N((e_-,\ c))& \ \ \
\ \ \ \  e_-= \sup \{x<c:\ x \in
E\}, \\
&r'=N((c,\ e_+))&\ \ \ \ \ \ \ e_+= \inf \{x>c:\ x \in E\}.
\end{eqnarray*}

The minimality of $D$ indicates $d_- \le e_- <c < e_+ \le d_+$,
which implies that $\ell \le \ell'$ and $r \le r'$.

So $R_D$ is the minimal renormalization.

(3) In order to describe the structure of $D$, we can consider the
following three cases, which cover all possible cases.
\begin{itemize}
\item Case A:\ \  $c \in D$,

\item Case B:\ \  $c \notin D$ and $D \cap (P_L,\ P_R)=\emptyset$,

\item Case C:\ \  $c \notin D$ and $D \cap (P_L,\ P_R)\neq
\emptyset$.

\end{itemize}

{\bf Case A:}\  If $c \in D$, the complete invariancy of $D$,
together with Lemma \ref{completely invariant}, implies $D=I$, which
is equivalent to $f$ is prime.

{\bf Case B:}\  If $c \notin D$ and $D\cap (P_L,\ P_R)=\emptyset$,
it follows from the proof of Case 2 that $D=O$. In this case, one
can check easily that $d_-=P_L$ and $d_+=P_R$ in the definition of
$R_D$. By Lemma \ref{equalities} we know that $N((P_L,\ c))=N((c,\
P_R))=\kappa$. It follows $R_D$ is periodic.

Conversely, assume that the minimal renormalization $R_D$ is
periodic. Follows from the definition of $R_D$, we know that the
renormalization interval of $R_D$ is $(f^{\kappa}(c+),\
f^{\kappa}(c-)) \subseteq (P_L,\ P_R)$ and the critcal interval of
$R_D$ is $(P_L,\ P_R)$. Consider the critical interval $(P_L,\
P_R)$, it follows $D\cap (P_L,\ P_R)=\emptyset$. So we get $D=O$ as
in Case 2 of the proof of (1).

{\bf Case C:}\  If $c \notin D$ and $D\cap (P_L,\ P_R)\neq
\emptyset$, it is necessary to prove $D$ is a Cantor set, the
equivalence between $D$ is a Cantor set and $R_D$ is not periodic is
obvious.

By Lemma \ref{completely invariant}, $c\notin D$ implies $D$ is
nowhere dense. Now we show that $D$ is perfect, i.e., $D=D'$.

Since $D=O$ is equivalent to
\begin{equation} \label{equality2}f^{\kappa}((P_L,\ P_R)) = (P_L,\
P_R),
\end{equation}
$D\cap (P_L,\ P_R)\neq \emptyset$ implies (\ref{equality2}) is not
true. Without loss of generality, we can suppose that
$f^{\kappa}(c+) < P_L$. Then there exists $y_1 \in (c,\ P_R)$ such
that $f^{\kappa}(y_1)=P_L$. And there exists $y_2 \in (c,\ P_R)$ and
$y_2>y_1$ such that $f^{\kappa}(y_2)=y_1$, i.e., $f^{2
\kappa}(y_2)=P_L$. Repeat the above arguments,  we can obtain an
increasing sequence $\{ y_n \}$ in $(c,\ P_R)$ such that $f^{n
\kappa}(y_n)=P_L$ and $y_n \to P_R$ as $n \to \infty$. Since
$\{y_n\}$ are preimages of $P_L$, we know that $\{y_n\} \subset D$.
It follows $P_R$ is a limit point of $D$, i.e. $P_R \in D'$. By
Lemma \ref{completely invariant} we know that $D'$ is backward
invariant, so $\cup_{n \ge 0} f^{-n}(P_R) \subset D'$. Therefore, $D
\subset D'$, $D$ is perfect.

The proof of Theorem B is complete. $\hfill \Box$

\vspace{0.5cm}

\subsection{The locally eventually onto property}

By Definition \ref{leo}, a Lorenz map is  {\it locally eventually
onto} ({\it l.e.o.}) if for any open interval $U$, there exists
positive integer $n$ depending on $U$, such that $\bigcup_{i =0}^n
f^i(U)=I$.

\begin{corollary}\label{leo-1}
Let $f$ be an expanding Lorenz map on $[a,\ b]$ with
$\kappa<\infty$. Then $f$ is {\it l.e.o.} if and only if it is
prime.
\end{corollary}

\begin{proof}
If $1<\kappa(f)<\infty$, by Lemma \ref{periodic orbit}, $f$ is {\it
l.e.o.} if and only if $D=I$. By Theorem B, $f$ is prime if and only
if $D=I$. So $f$ is {\it l.e.o.} is equivalent to it is prime.

If $\kappa(f)=1$, it is easy to see that $f$ is prime, and $f$ is
{\it l.e.o.}, because either $a$ or $b$ is fixed.
\end{proof}

Glendinning and Sparrow  \cite{GS} described the {\it locally
eventually onto } {(\it l.e.o.)} property as follows: $f$ is said to
be {\it l.elo.} if for each open interval $U$, there exists
subintervals $U_1$, $U_2$ of $U$, and positive integers $n_1$, $n_2$
such that $f^{n_1}$ and $f^{n_2}$ map $U_1$ and $U_2$
homeomorphically to $(a,\ c)$ and $(c,\ b)$, respectively.  The
following proposition relates two  definitions of ${\it l.e.o.}$.

\begin{proposition}\label{equivalence}
Our definition of {\it l.e.o.} coincides with which of Glendinning
and Sparrow in \cite{GS} when $\kappa\le 2$.
\end{proposition}

\begin{proof}
It is necessary to show that our definition of {\it l.e.o.} reduces
to Glendinning and Sparrow's definition when $\kappa \le 2$. The
converse is trivial.

Now suppose $f$ is prime and $\kappa \le 2$. There are two cases:
$\kappa=1$ and $\kappa=2$.

If $\kappa=1$, at least one of the following holds:
$$f(a)=a\ \ and\ \  f(b)=b.$$
Without loss of generality, we suppose $f(a)=a$. For any open
interval $U=(x,\ y)$, let $z_0 $ be the point in $U$ such that
$f^{N(U)}(z_0)=c$, and $z_1$ be the point in $(z_0,\ y) \subset U$
such that $f^{N((z_0,\ y))}(z_1)=c$. By the definition of Lorenz map
and $f(a)=a$ we obtain
$$f^{N((z_0,\ y))+1}((z_0,\ z_1))=(a,\ b).$$
So there exists positive integers $n$ and a subinterval $V \subseteq
U$ such that $f^n$ maps $V$ to $(a,\ b)$ homeomorphically, which
implies that $f$ is locally eventually onto.

For the case $\kappa=2$. Suppose $f$ is prime, let $P_L<c<P_R$ be
the $2-$periodic points.  By Lemma \ref{equalities}, we know that
$N((P_L,\ c))=N((c,\ P_R))=2$, so $N((P_R,\ b))=1$ because $f((P_L,\
c))=(P_R,\ b)$. Let $x_1$ be the point in $(P_R,\ b)$ such that
$f^2(x_1)=c$, $y_1$ be the point in $(P_R,\ b)$ such that
$f(y_1)=c$. Consider the interval $J_1=(x_1,\ y_1)$, one can check
that $f^2(J_1)=(c,\ b) \supset J_1$. There exists an subinterval
$J_2 \subset J_1$ so that $f^2(J_2)=J_1$. So we can obtain a
sequence of nested intervals $\{J_n\}_n$, $J_n=(x_n,\ y_n)$ satisfy:
$$J_{n+1} \subset J_n,\ \ \ f^2(J_{n+1})=J_n, \ \ \ f^{2n}(J_n)=(c,\ b),\ \  n=1,\ 2,\, \cdots.$$
Since $\{x_n\}$ and $\{y_n\}$ are monotone and $f$ is expanding, the
length of $|J_n| \to 0$ as $n \to \infty$.

Now we prove that $f$ is {\it l.e.o.} in the sense of Glendinning
and Sparrow. It is necessary to check the {\it l.e.o.} conditions
for intervals containing $P_R$, because $f$ is prime implies that
any open interval contains a subinterval which can be mapped
homeomorphically to an open interval containing $P_R$. For any open
interval $F$ containing $P_R$, we can find subinterval $J_i \subset
F$, which implies that $f^{2i}$ maps $J_i$ homeomorphically to $(c,\
b)$ by the construction of $\{J_n\}$. Furthermore, $(c,\ b)$
contains an interval $(c,\ y_1)$, which can be mapped by $f$
homeomorphically to $(a,\ c)$. So $J_i$ contains a subinterval
$(x_i, z_i)$ such that $f^{2i+1}((x_i,\ z_i))=(a,\ c)$. Hence, $f$
is {\it l.e.o.} in the sense of Glendinning and Sparrow.
\end{proof}

\begin{remark} The exact formulation of {\it l.e.o.} varies in the
literatures. For the definition we use, we mention the following:
\begin{enumerate}

\item The {\it l.e.o.} property is just the strongly transitive
property in Parry \cite{P3}.

\item It agrees with the one in \cite{GS}, when $\kappa \le 2$;

\item $f$ is prime if and only if $f$ is {\it l.e.o.};

\item The {\it l.e.o.} property of expanding Lorenz map comes from
the {\it l.e.o.} property of the periodic orbit with minimal period.
According to Lemma \ref{periodic orbit},\ for expanding Lorenz map
$f$, the minimal completely invariant closed set $D$ of $f$ admits
the {\it l.e.o.} property: for each open interval $U$ satisfying $U
\cap D \neq \emptyset$, there exists integer $n>0$ so that
$\bigcup_{i=0}^{n}f^i(U)=I$. As a result, $f$ is {\it l.e.o.} if and
only if $D=I$.

\end{enumerate}

\end{remark}

\section{Consecutive renormalizations: characterization of $\alpha$-limit set}  \ \

Thanks to Theorem B, the minimal renormalizaion of renormalizable
expanding Lorenz map always exists. We can define a renormalization
operator $R$ from the set of renormalizable expanding Lorenz maps to
the set of expanding Lorenz maps. For each renormalizable expanding
Lorenz map, $Rf:=R_Df$, where $D$ is the minimal proper completely
invariant closed set of $f$. Obviously, $Rf$ is also expanding. If
$Rf$ is renormalizable, we can obtain $R^2f:=R(Rf)$. In this way, we
define $R^nf$ as the minimal renormalization of $R^{n-1}f$ if
$R^{n-1}f$ is renormalizable. If the renormalization process can
proceed $m$ times, we say that $f$ is $m$ ($0 \le m \le \infty$)
{\it times renormalizable}. If $f$ is $m$-renormalizable, then
$\{R^i f\}_{i=1}^m$ are all the renormalizations of $f$. We call
$R^if$ the $i$th renormalization of $f$. The process of consecutive
renormalizations can be used to characterize all the $\alpha$-limit
sets and nonwandering set of expanding Lorenz map.

\subsection{$\alpha$-limit set}

\begin{lemma}\label{alpha-limit set}
Let $f$ be an expanding Lorenz map. Each proper completely invariant
closed set of $f$ is an $\alpha$-limit set.
\end{lemma}

\begin{proof} Suppose $f$ is $m$-renormalizable ($0 \le m \le \infty$),
with renormalization intervals $[a_i,\ b_i]$, $i= 1,\ \cdots, m$.
There are $m$ proper completely invariant closed sets for $f$,
\begin{equation}\label{alpha-limit set1}
E_i=\{x: orb(x) \cap (a_i,\ b_i)=\emptyset \},  \ \ \ \ i= 1,\
\cdots, m.
\end{equation}
We have $$E_1 \subset E_{2}\subset \cdots \subset E_m
$$ because $[a_i,\ b_i] \supset [a_{i+1},\ b_{i+1}]$, $0 <i < m$.

Now we prove that $E_i$ is an $\alpha$-limit set of $f$ for $0 <i
\le m$. Put $e_-^i=\sup \{x \in E_i: x<c\}$. According to Theorem A
we know that $e_-^i$ is periodic. By Lemma \ref{completely
invariant}, \ $\alpha(e_-^i)$ is indeed a completely invariant
closed set, and \ $e_-^i \in \alpha(e_-^i)$. We must have
$\alpha(e_-^i)=E_k$ for some $k=1, 2, \ldots, m$, because $f$ admits
exact $m$ proper completely invariant closed sets.

Since
$$
(e^{i-1}_-,\ e^{i-1}_+) \supset (a_{i-1},\ b_{i-1}) \supset
(e^{i}_-,\ e^{i}_+) \supset (a_{i},\ b_{i})  \supset (e^{i+1}_-,\
e^{i+1}_+) \supset (a_{i+1},\ b_{i+1}),
$$
by the definition of $E_i$ and $E_{i+1}$, we know that $e^{i}_-
\notin E_{i-1}$ and $e_-^{i+1} \in E_{i+1}\backslash E_i$.

Observe that $e^i_- \in \alpha(e^i_-)$ and  $e^{i}_- \notin E_{i-1}$
indicate that $k \ge i$, and $e_-^{i+1} \in E_{i+1}\backslash E_i$
implies $k<i+1$, we conclude that $k=i$, $i.e.$,
$\alpha(e_-^i)=E_i$. Hence, $E_i$ is an $\alpha$-limit set.
\end{proof}

\vspace{0.5cm}



\subsection{Proof of Theorem C}


(1), By Lemma \ref{completely invariant} we know that each
$\alpha$-limit set is completely invariant. And by Lemma
\ref{alpha-limit set} each completely invariant set is an
$\alpha$-limit set. So completely invariant closed set and
$\alpha$-limit set of $f$ are the same thing in different names. If
$f$ is $m$-renormalizable, then $f$ has exact $m$ proper
$\alpha$-limit sets. Follows from the proof of Lemma
\ref{alpha-limit set}, all the $\alpha$-limit sets are
$\{E_i\}_{i=1}^m$ defined in (\ref{alpha-limit set1}), and
$$
\emptyset=E_0\subset E_1 \subset E_2  \subset \cdots \subset E_m
\subset I.
$$

(2), At first we prove that if the $i$th $(0 <i \le m< \infty)$
renormalization is periodic, then $E_{i}'=E_{i-1}$.

Suppose $g=R^{i-1}f$. $g$ is an expanding Lorenz map on $[a_{i-1},\
b_{i-1}]$ with discontinuity $c$. Denote $\kappa_1$ as the minimal
period of periodic points of $g$, $O_1$ as the $\kappa_1$-periodic
orbit of $g$, and $P_L'$ and $P_R'$ are two adjacent
$\kappa_1$-periodic point of $g$ with $P_L'<c<P_R'$. By Lemma
\ref{completely invariant} and the proof of Lemma \ref{alpha-limit
set} we know that $E_i=\overline{\bigcup_{n \ge 0}f^{-n}(P_L')}$.

Put $$e^{i-1}_-=\sup\{x \in E_{i-1}, \ x <c\},\ \ \ \ \
e^{i-1}_+=\inf\{x \in E_{i-1}, \ x >c\},$$
$$
\ell'=N((e^{i-1}_-,\ c)),\ \ \ \ \ \ \ r'=N((c,\ e^{i-1}_+)).
$$
According to the definition of minimal renormalization, we have
$e^{i-1}_-<a_{i-1} \le P_L'<c<P_R'\le b_{i-1}<e^{i-1}_+$.

By assumption, the minimal renormalization of $g$ is periodic, it
follows from Theorem B that the minimal completely invariant closed
set of $g$ is $O_1$, which implies $P_L'$ is an isolated point of
$E_i$. So $E_i' \neq E_i$.

Observe that $f^{\ell'}((e^{i-1}_-,\ c))=(e^{i-1}_-,\ b_{i-1})$,
there exists a decreasing sequence $\{x_n\}$ in $E_{i-1}\cap
(e^{i-1}_-,\ c)$ such that
$$
f^{\ell'}(x_{1})=P_R',\ \ \ \ f^{\ell'}(x_{n+1})=x_n,\ \ \ \ n=1,\
2,\ \cdots
$$
and $x_n \to e^{i-1}_-$ as $n \to \infty$. So $e^{i-1}_- \in E_i'$.

By Lemma \ref{completely invariant} we know $E_i'$ is also a
completely invariant closed set, we have
$$E_{i-1}=\overline{\bigcup_{n \ge 0}f^{-n}(e^{i-1}_-)} \subseteq E_i' \neq E_{i}.$$
It follows $E_i'=E_{i-1}$.

Now we show that if the $i$th renormalization $R^if$ is not
periodic, then $E_i$ a Cantor set. From the proof of first part, we
know that $E_i=\overline{\bigcup_{n \ge 0}f^{-n}(P_L')}$. Since the
$i$th renormalization is not periodic, the minimal completely
invariant closed set of $R^{i-1}f$ is a Cantor set. So $E_i$ admits
no isolated point in $[a_{i-1},\ b_{i-1}]$. $E_i' \cap [a_{i-1},\
b_{i-1}] \neq \emptyset$, which implies that $E_i'=E_i$.

(3) Now we are ready to characterize the $\alpha$-limit set of every
point in $I$. At first, we describe the set $\{x \in I,
\alpha(x)=D\}$, where $D$ is the minimal completely invariant closed
set of $f$.

{\bf Claim: }$\alpha(x)=D$ if and only if $x \notin orb([a_1,\
b_1])$, where $[a_1,\ b_1]$ is the renormalization interval of the
minimal renormalization $R_D$.

Suppose the minimal renormalization $Rf=R_Df:=(f^{\ell}, f^r)$. It
follows that
$$
orb([a_1,\ b_1])=\bigcup_{n \ge 0}f^n([a_1,\
b_1])=\left(\bigcup_{n=0}^{\ell-1}f^n([a_1,\ c])\right) \bigcup
\left(\bigcup_{n=0}^{r-1}f^n([c,\ b_1])\right)
$$
is the union of finite closed intervals, and $orb([a_1,\ b_1])$ is
forward invariant under $f$.

Since $D$ is the minimal completely invariant closed of $f$, by
Lemma \ref{alpha-limit set}, $D$ is also the minimal $\alpha$-limit
set of $f$. So $\alpha(x) \supset D$ for all $x \in I$.

Let $D_1$ be the minimal completely invariant closed set of the
minimal renormalization $R_Df$. It follows that  $D_1 \cap
D=\emptyset$, and $D_1 \subset E_2$. If $x\notin orb([a_1,\ b_1])$,
then $f^{-n}(x) \cap orb([a_1,\ b_1])=\emptyset$ because $orb([a_1,\
b_1])$ is forward invariant under $f$. So $\alpha(x)$ is disjoint
with the interior of $orb([a_1,\ b_1])$, which indicates $\alpha(x)
\cap D_1=\emptyset$. Hence, $\alpha(x) \neq E_2$, i.e.,
$\alpha(x)=D=E_1$.

On the other hand, by the minimality of $D_1$, $\alpha(x, R_Df)=D_1$
for all $x \in [a_1,\ b_1]$. For $x \in [a_1,\ b_1]$, since $orb(x,
R_Df)=orb(x, f) \cap [a_1,\ b_1]$, we see that $\alpha(x) \supset
\alpha(x, R_Df)$. So $\alpha(x) \cap D_1 \neq \emptyset$, which
implies that $\alpha(x) \neq D$ for $x \in [a_1,\ b_1]$. Notice that
$\alpha(x) \subseteq \alpha(f(x))$, we conclude $\alpha(x) \neq D$
for all $x \in orb([a_1,\ b_1])$.

The proof of the Claim is complete.

For $0 \le i\le m$, we denote $[a_i,\ b_i]$ as the renormalization
interval of the $i$th renormalization $R^if$, and $D_i$ as the
minimal completely invariant closed set of $R^if$.

By the Claim we know that $\alpha(x)=E_1$ if and only if $$x \in I
\backslash orb([a_1,\ b_1])= orb([a_0,\ b_0])\backslash orb([a_1,\
b_1]).$$

For the case $i=2 \le m$, we consider the map $Rf:=R_Df$ on $[a_1,\
b_1]$. According to the Claim, we obtain that $\alpha(x,\ Rf)=D_1$
if and only if $x\notin orb([a_2,\ b_2]$. It follows that
$\alpha(x)=E_2$ if and only if
$$
x\in orb([a_1,\ b_1])\backslash orb([a_2,\ b_2]).
$$

Repeat the above arguments, we conclude $\alpha(x)=E_i$ if and only
if
$$x \in orb([a_{i-1},\ b_{i-1}])\backslash orb([a_i,\ b_i]) \ \  for
\ \ 0<i\le m.$$

If $m<\infty$, $R^mf$ is prime on $[a_m,\ b_m]$, $\alpha(x,
R^mf)=[a_m,\ b_m]$ for all $x\in [a_m,\ b_m]$. By Lemma
\ref{completely invariant}, the completely invariant closed set
containing $[a_m,\ b_m]\ni c$ is $I$, we conclude that $\alpha(x)=I$
for all $x\in orb([a_m,\ b_m])$.

For  the case $m=\infty$, put $A=\cap_{i\ge 1}^m orb([a_i,\ b_i])$,
it is known that $A:=\overline{orb(c+)}=\overline{orb(c-)}$ (cf.
\cite{GS} \cite{KM}), which is a Cantor set. Since $c\in A$, the
completely invariant closed set containing $A$ is $I$. As a result,
$\alpha(x)=I$ for all $x\in A$. $\hfill \Box$


\vspace{0.5cm}

\subsubsection{Example: $\alpha$-limit set with given depth} \ \

We can use Theorem C to construct countable $\alpha$-limit set with
given depth.

Consider the piecewise linear symmetric Lorenz map: $1<a \le 2$,
\begin{equation}\label{LinearLorenzMap-0}
f_a(x)= \left \{ \begin{array}{ll}
ax+1-\frac{1}{2}a & x \in [0,  \frac{1}{2}) \\
a(x-\frac{1}{2}) & x \in (\frac{1}{2},  1].
\end{array}
\right.
\end{equation}

According to Glendinning \cite{G} and Palmer \cite{Pal} , $f_a$ can
only be periodically renormalized finite times. Suppose $a \in
(2^{2^{-(m+1)}},\ 2^{2^{-m}}]$, Parry \cite{P4} proved that $f_a$
can be (periodically) renormalized $m$ times. In this case, by
Theorem A and Theorem C, $f_a$ has exact $m$ different
$\alpha$-limit sets. Let $p_i$ be one of the $2^i$-periodic point of
$f_a$,
$$
E_i=\overline{\bigcup_{n \ge 0}f_a^{-n}(p_i)}, \ \ \ \ \ i=1,\
\cdots, m.
$$
Then $\{E_i\}_{i=1}^m$ is the cluster of $\alpha$-limit sets of
$f_a$. Moreover, according to Theorem C, $E_n^{(i)}=E_{n-i}$, $i=1,\
2,\ \ldots, \ m$. So $E_m$ is a countable closed set and the $m$the
derived set $E_m^{(m)}$ is empty. The depth of $E_m$ is $m$.




\vspace{0.2cm}

\subsection{Nonwandering set} \ \

The following Lemma \ref{mminimal} indicates that the dynamics on
the minimal completely invariant closed $D$ is indecomposable.

\begin{lemma}\label{mminimal}
Let $f$ be an expanding Lorenz map with $1 \le \kappa< \infty$, $D$
be its minimal completely invariant closed set. Then $f: D \to D$ is
{\it l.e.o.}, and $\Omega(f|_D)=D$.
\end{lemma}

\begin{proof}
By Theorem B, there are three cases: $D=O$, $D=I$ and $O \subset D
\subset I$, where $O$ is the unique $\kappa$-periodic orbit. If
$D=O$ or $D=I$, applying Theorem B, it is easy to see $f: D \to D$
is {\it l.e.o.}, and $\Omega(f|_D)=D$.

For the case $O \subset D \subset I$, we know that $D$ is  a
completely invariant Cantor set of $f$. We shall prove that $f: D
\to D$ is {\it l.e.o.}. Suppose $A$ is an open set of $D$ (in the
induced topology from $I$), there exists an open set $U$ of $I$ such
that $A=U \cap D$. By the {\it l.e.o.} property of $O$, there is
positive integer $N$ such that $\bigcup_{n=0}^N f^n(U)=I$. It
follows from Lemma \ref{completely invariant} that $f^n(U \cap
D)=f^n(U) \cap D$. We have
$$
\bigcup_{n=0}^n f^n(A)=\bigcup_{n=0}^n f^n(D \cap
U)=(\bigcup_{n=0}^n f^n(U)) \cap D=D,
$$
which implies that $f: D \to D$ is {\it l.e.o.}. As a result,
$\Omega(f|_D)=D$.
\end{proof}

 \vspace{0.3cm}
Now we prove the nonwandering set decomposition of expanding Lorenz
map. As mentioned before, Glendinning and Sparrow \cite{GS} gave a
decomposition based on kneading theory.

 \vspace{0.3cm}

\begin{proposition} \label{Proposition D}  Let $f$ be an $m$-renormalizable
($0 \le m \le \infty$) expanding Lorenz map and
$$\emptyset=E_0 \subset E_1 \subset E_2 \subset \cdot \cdot \cdot \subset E_m \subset I$$
be all the $\alpha$-limit sets of $f$, $I_i=[a_i,\ b_i]$ be the
renormalization interval of the $i$th renormalization $R^if$, and
$D_i$ is the minimal completely invariant closed set of $R^if$.

Then there is a canonical decomposition of the nonwandring set
$\Omega(f)$ of $f$ into $m$-invariant closed set $\Omega_i$
($i=1,\cdots, m$) and an attractor $A$
\begin{equation}\label{nonwandering set}
\Omega(f)=\bigcup_{i=1}^m\Omega_i   \cup A.
\end{equation}

This decomposition has the following properties:

\begin{enumerate}

\item $\Omega_i:=E_i \cap orb(I_{i-1})=orb(D_{i-1})$, $1\le i
\le m$ and $f|_{\Omega_i}$ is {\it l.e.o.}.  $\Omega_i$ is either a
periodic orbit or a Cantor set depending on wether the
renormalization $R^if$ is periodic or not.

\item $A$ is the attractor of $f$: $\omega(x) \subseteq A$ for $x \notin E_{\infty}:=\bigcup_{i\ge
0}E_i$. $f|_A$ is {\it l.e.o.}. Moreover, $A=\bigcap_{i=0}^m
orb([a_i,\ b_i])$:\ $A=I$ if $m=0$, $A$ is a finite union of closed
intervals if $0<m<\infty$, and $A$ is a Cantor set if $m=\infty$. In
the last case, $\omega(x)=A$ for $x \notin E_{\infty}$.
\end{enumerate}
\end{proposition}

\vspace{0.3cm}

Based on their renormalization theory on kneading invariants,
Glendinning and Sparrow \cite{GS} obtained the nonwandering set
decomposition (\ref{nonwandering set}). Our proof of the
decomposition is independent of kneading theory. We obtain the exact
expression of $\Omega_i$, and emphasize that $\Omega_i$ is
indecomposable: $f|_{\Omega_i}$ is {\it l.e.o.}.

\vspace{0.3cm}

\begin{proof}
If $m=0$,  $f$ is prime. By Theorem B and Theorem C, we know that
$f$ is ${l.e.o.}$, $A=I=\Omega(f)$, and $\alpha(x)=I$, $\forall x
\in I$.

Now suppose $m>0$, i.e., $f$ is renormalizable. By Theorem A and
Theorem C, all the completely invariant closed sets of $f$ are:
$$\emptyset=E_0 \subset E_1 \subset E_2 \subset \cdot \cdot \cdot \subset E_m \subset E_{m+1}=I.$$
($E_{m+1}=I$ is just a notation when $m=\infty$). We can decompose
$I=E_{m+1}$ as follows:
$$
I=(E_1\backslash E_0) \cup (E_2\backslash E_1) \cup \cdots \cup
(E_m\backslash E_{m-1}) \cup  (E_{m+1}\backslash E_m).
$$
Since $E_i$ is completely invariant, $E_i\backslash E_{i-1}$
($i=1,\ldots, m$) and $I\backslash E_m$ are invariant under $f$. It
follows that
\begin{eqnarray*}
\Omega(f)=\Omega(f) \bigcap \bigcup_{i=1}^{m+1}(E_i\backslash
E_{i-1})=\bigcup_{i=1}^{m+1}(\Omega(f) \cap (E_i\backslash
E_{i-1})):=\bigcup_{i=1}^m \Omega_i \cup A
\end{eqnarray*}
where $\Omega_i=\Omega(f) \cap (E_i\backslash E_{i-1})$ and
$A=\Omega(f) \cap (I\backslash E_{m})$.

For $0<i \le m$, the $(i-1)$th renormalization of $f$ is
\begin{equation}\label{ith renormal}
R^{i-1}f(x)= \left \{ \begin{array}{ll}
f^{\ell_{i-1}}(x) & x \in [a_{i-1},  c) \\
f^{r_{i-1}}(x) & x \in (c,  b_{i-1}].
\end{array}
\right.
\end{equation}
$D_{i-1}$ is the minimal completely invariant closed set of
$R^{i-1}f$.

In what follows, we only show that $$\Omega_i=orb(D_{i-1})=E_i \cap
orb ([a_{i-1},\ b_{i-1}])$$ in three steps. By Lemma \ref{mminimal}
we know that $f|_{\Omega_i}$ is {\it l.e.o.}. For the proof of
remain parts, see \cite{GS} or \cite{KM}.

{\bf Step 1:} $orb(D_{i-1})=E_{i} \cap orb([a_{i-1}, b_{i-1}]).$

By the definitions of $E_i,\ E_{i-1}$ and $D_{i-1}$,  $D_{i-1}
\subseteq E_i \cap [a_{i-1}, b_{i-1}]$. On the other hand, $x \in
[a_{i-1},\ b_{i-1}] \backslash E_i$ indicates $orb(R^{i-1}f, x) \cap
[a_i,\ b_i] \neq \emptyset$. By Lemma \ref{complimentary}, $x \notin
D_{i-1}$. We obtain $[a_{i-1},\ b_{i-1}]\backslash D_{i-1} \subseteq
[a_{i-1},\ b_{i-1}]$, which implies $D_{i-1} \supseteq E_i \cap
[a_{i-1},\ b_{i-1}]$. By Lemma \ref{completely invariant}, we get
the desired equality.

{\bf Step 2:} $orb(D_{i-1})\subseteq \Omega_i:=\Omega(f) \cap
(E_i\backslash E_{i-1}).$

By definitions, $orb(D_{i-1}) \subset E_i\backslash E_{i-1}$.

By Lemma \ref{mminimal}, we know that $R^{i-1}f|_{D_{i-1}}$ is {\it
l.e.o.}, and $\Omega(R^{i-1}f|_{D_{i-1}})=D_{i-1}$. Since $R^{i-1}f$
is the first return map of $f$ on the renormalization interval
$I_{i-1}:=[a_{i-1},\ b_{i-1}]$, we have
$$
orb(x, R^{i-1}f)=orb(x, f) \cap I_{i-1},\ \ \ \ \forall x \in
I_{i-1}.
$$
It follows that $D_{i-1} \subset \Omega(R^{i-1}f) \subset
\Omega(f)$, and $orb(D_{i-1}) \subset \Omega(f)$ because $\Omega(f)$
is invariant under $f$.

{\bf Step 3:}  $orb ([a_{i-1},\ b_{i-1}]) \supseteq
\Omega_i:=\Omega(f) \cap (E_i\backslash E_{i-1}).$

It is necessary to show that any point in $E_i\backslash
(E_{i-1}\cup orb ([a_{i-1},\ b_{i-1}])$ is wandering. Suppose $x \in
E_i$, and $x \notin (E_{i-1}\cup orb ([a_{i-1},\ b_{i-1}])$. $x
\notin E_{i-1}$ implies the orbit of $x$ will go into $(a_{i-1},\
b_{i-1})$, and stay in the forward invariant closed set
$orb([a_{i-1},\ b_{i-1}])$ forever. So $x$ is wandering because $x
\notin orb([a_{i-1},\ b_{i-1}])$.
\end{proof}

\vspace{0.3cm}

{\bf Appendix: Periodic renormalization} \ \

\vspace{0.2cm}

We collect some facts for periodic renormalization in this appendix.


Let $f$ be an expanding Lorenz map on $[a,\ b]$, $\kappa$ is the
minimal period of $f$, $O$ is the unique $\kappa$-periodic orbit,
$P_L$ is the largest $\kappa-$periodic point in $[a, c)$ and $P_R$
be the smallest $\kappa-$periodic point in $(c, b]$, $D$ is the
minimal completely invariant closed set of $f$. If $f$ is
$m-$renormalizable, $0 \le i \le m$, $R^i$ if is the $i$th
renormalization of $f$, and $E_i$ is the completely invariant closed
set corresponds to $R^if$, $\Omega_i=E_i \cap \Omega(f)$.

\begin{enumerate}

\item The minimal renormalization $Rf$ is periodic if and only if
$D=O$ (Theorem B).

\item  $Rf$ is periodic if and only if
\begin{equation} \label{periodic renormalization}
[f^{\kappa}(c+),\ f^{\kappa}(c-)] \subseteq [P_L,\ P_R].
\end{equation}

\item One can check if $Rf$ is periodic or not in following
steps:
\begin{itemize}
\item Find the minimal period $\kappa$ of $f$ by considering the preimages of
$c$, see Lemma \ref{minimal period};

\item Find the $\kappa$-periodic orbit;

\item Check if the inclusion (\ref{periodic renormalization}) holds
or not.

\end{itemize}

\item  $Rf$ is periodic if and only if
the rotational interval of $f$ is degenerated to a rational point
\cite{AD, M}.

\item $R^if$ is periodic if and only if $E_i$ admits isolated
points.

\item $R^if$ is periodic if and only if $\Omega_i$ is consists of a
periodic orbit.

\item $R^if$ is periodic if and only if the topological entropy  $h(f|_{\Omega_i})=0$.

\item If the first $k (k\le m)$ renormalizations are all periodic,
then $E_i=E_k^{(k-i)}$, $i=0,1,\ldots, k$. The depth of $E_i$ is
$i$.

\item Suppose $f$ is a $\beta$-transformation. $f$ can only be renormalized periodically, $i.e.$,  each renormalization of $f$ is
periodic. Since $\beta$-transformation is finitely renormalizable,
one can obtain all of the renormalizations of $f$ in finite steps.

\item An expanding Lorenz map $f$ is conjugated to a
$\beta$-transformation if and only if $f$ is finitely renormalizable
and each renormalization of $f$ is periodic \cite{G, Pal}.

\item A piecewise linear Lorenz map that expand on average is
conjugate to $\beta$-transformation \cite{CD2}.

\item For $1<a \le 2$, put
\begin{equation}\label{LinearLorenzMap-0}
f_a(x)= \left \{ \begin{array}{ll}
ax+1-\frac{1}{2}a & x \in [0,  \frac{1}{2}) \\
a(x-\frac{1}{2}) & x \in (\frac{1}{2},  1].
\end{array}
\right.
\end{equation}
If  $a \in (2^{2^{-(m+1)}},\ 2^{2^{-m}}]$, then $f_a$ is
$m$-renormalizable \cite{P4}.
\end{enumerate}

\vspace{0.3cm}

{\bf Acknowledgements:}  This work is partially supported by a grant
from the Spanish Ministry (No. SB2004-0149) and  grants from NSFC
(Nos. 60534080, 70571079) in China. The author thanks Centre de
Recerca Matem$\grave{a}$tica for the hospitality and facilities. The
author thanks Prof. Lluis Alsed$\grave{a}$ and Dr. Hongfei Cui for
many discussions and suggestions.

\end{document}